\documentclass[12pt]{amsart}

\usepackage{amsmath}
\usepackage{amsxtra}
\usepackage{amscd}
\usepackage{amsthm}
\usepackage{amsfonts}
\usepackage{amssymb}
\usepackage{eucal}
\usepackage[matrix,arrow,curve]{xy}
\usepackage[dvips]{graphicx}
\usepackage[dvips]{graphics}
\usepackage[T2A]{fontenc}
\usepackage[cp866]{inputenc}

\newcommand{\Res}{\mathop{\sf Res}\nolimits}

\theoremstyle{plain}
\newtheorem{Thm}[subsection]{Theorem}
\newtheorem{Cor}[subsection]{Corollary}
\newtheorem{Lem}[subsection]{Lemma}
\newtheorem{Prop}[subsection]{Proposition}
\newtheorem{Conj}[subsection]{Conjecture}
\newtheorem{Ex}[subsection]{Example}

\theoremstyle{definition}
\newtheorem{Def}[subsection]{Definition}

\theoremstyle{remark}

\newtheorem{Rem}[subsection]{Remark}

\newtheorem{conj}{Conjecture}

\numberwithin{equation}{section}

\newif\ifShowLabels
\ShowLabelstrue
\newdimen\theight
\def\TeXref#1{%
    \leavevmode\vadjust{\setbox0=\hbox{{\tt
        \quad\quad  {\small \rm #1}}}%
    \theight=\ht0
    \advance\theight by \lineskip
    \kern -\theight \vbox to
    \theight{\rightline{\rlap{\box0}}%
    \vss}%
    }}%

\ShowLabelsfalse

\renewcommand{\sec}[2]{\section{#2}\label{S:#1}%
    \ifShowLabels \TeXref{{S:#1}} \fi}
\newcommand{\ssec}[2]{\subsection{#2}\label{SS:#1}%
    \ifShowLabels \TeXref{{SS:#1}} \fi}

\newcommand{\refs}[1]{Section ~\ref{S:#1}}

\newcommand{\reft}[1]{Theorem ~\ref{T:#1}}
\newcommand{\refl}[1]{Lemma ~\ref{L:#1}}
\newcommand{\refp}[1]{Proposition ~\ref{P:#1}}

\newcommand{\refr}[1]{Remark ~\ref{R:#1}}
\newcommand{\refe}[1]{\eqref{E:#1}}

\newenvironment{thm}[1]%
    { \begin{Thm} \label{T:#1}  \ifShowLabels \TeXref{T:#1} \fi }%
    { \end{Thm} }

\renewcommand{\th}[1]{\begin{thm}{#1} \sl }
\renewcommand{\eth}{\end{thm} }

\newenvironment{lemma}[1]%
    { \begin{Lem} \label{L:#1}  \ifShowLabels \TeXref{L:#1} \fi }%
    { \end{Lem} }
\newcommand{\lem}[1]{\begin{lemma}{#1} \sl}
\newcommand{\elem}{\end{lemma}}

\newenvironment{propos}[1]%
    { \begin{Prop} \label{P:#1}  \ifShowLabels \TeXref{P:#1} \fi }%
    { \end{Prop} }
\newcommand{\prop}[1]{\begin{propos}{#1}\sl }
\newcommand{\eprop}{\end{propos}}

\newenvironment{corol}[1]%
    { \begin{Cor} \label{C:#1}  \ifShowLabels \TeXref{C:#1} \fi }%
    { \end{Cor} }
\newcommand{\cor}[1]{\begin{corol}{#1} \sl }
\newcommand{\ecor}{\end{corol}}

\newenvironment{defeni}[1]%
    { \begin{Def} \label{D:#1}  \ifShowLabels \TeXref{D:#1} \fi }%
    { \end{Def} }
\newcommand{\defe}[1]{\begin{defeni}{#1} \sl }
\newcommand{\edefe}{\end{defeni}}

\newenvironment{remark}[1]%
    { \begin{Rem} \label{R:#1}  \ifShowLabels \TeXref{R:#1} \fi }%
    { \end{Rem} }
\newcommand{\rem}[1]{\begin{remark}{#1}}
\newcommand{\erem}{\end{remark}}

\newenvironment{conjec}[1]%
    { \begin{Conj} \label{Co:#1}  \ifShowLabels \TeXref{Co:#1} \fi }%
    { \end{Conj} }
\renewcommand{\conj}[1]{\begin{conjec}{#1} \sl }
\newcommand{\econj}{\end{conjec}}

\newenvironment{example}[1]%
    { \begin{Ex} \label{Exx:#1}  \ifShowLabels \TeXref{Exx:#1} \fi }%
    { \end{Ex} }
\newcommand{\ex}[1]{\begin{example}{#1} \sl }
\newcommand{\eex}{\end{example}}

\newcommand{\eq}[1]%
    { \ifShowLabels \TeXref{E:#1} \fi
       \begin{equation} \label{E:#1} }
\newcommand{\eeq}{ \end{equation} }

\newcommand{\prf}{ \begin{proof} }
\newcommand{\epr}{ \end{proof} }

     \newcommand\Gam{\Gamma}

\newcommand\zet{\zeta}
     
\newcommand\iot{\iota}

\newcommand\lam{\lambda}        \newcommand\Lam{\Lambda}
\newcommand\sig{\sigma}     

\newcommand\ome{\omega}     \newcommand\Ome{\Omega}

\newcommand\calW{{\mathcal{W}}}

\newcommand\bfk{{\mathbf k}}

\newcommand\bft{{\mathbf t}}

\newcommand\bfx{{\mathbf x}}        
\newcommand\bfy{{\mathbf y}}        
\newcommand\bfz{{\mathbf z}}        

\newcommand\CC{\mathbb{C}}

 \newcommand\grg{{\mathfrak{g}}}

\newcommand\sdp{\times \hskip -0.3em {\raise 0.3ex
\hbox{$\scriptscriptstyle |$}}} 

\newcommand\Gr{\operatorname{Gr}}

\newcommand\Id{\operatorname{Id}}

\newcommand\Int{\operatorname{Int}}

\newcommand\sgn{\operatorname{sgn}}

\newcommand\tilE{{\widetilde{E}}}
\newcommand\tilf{{\widetilde{f}}}
\newcommand\tilF{{\widetilde{F}}}

\newcommand\tilG{{\widetilde{G}}}

\newcommand\x{\times}
\newcommand\ten{\otimes}

\renewcommand{\Id}{\text{Id}}

\newcommand\nc{\newcommand}

\newcommand{\iso}{{\stackrel{\sim}{\longrightarrow}}}

\nc\aff{\operatorname{aff}}
\nc\oGr{\overline{\Gr}}
\nc\Bun{\operatorname{Bun}}
\nc\hgrg{\widehat{\grg}}
\renewcommand\Int{\operatorname{Int}}
\nc\bInt{\overline{\Int}}
\nc\hatLam{\widehat{\Lam}}
\nc\bmu{\overline{\mu}}
\nc\bnu{\overline{\nu}}
\nc\blambda{\overline{\lam}}

\nc\ocalW{\overline{\calW}}
\nc\pos{\operatorname{pos}}
\nc\IH{\operatorname{IH}}
\nc\Rep{\operatorname{Rep}}
\nc\Gal{\operatorname{Gal}}
\nc{\tilGr}{\widetilde{\Gr}}

\nc\Pic{\operatorname{Pic}}
\nc\pa{\partial}
\nc\Na{\nabla}

\nc{\HC}{{\mathcal{HC}}}
\nc{\on}{\operatorname}
\nc{\BA}{{\mathbb{A}}}
\nc{\BC}{{\mathbb{C}}}
\nc{\BG}{{\mathbb{G}}}
\nc{\BM}{{\mathbb{M}}}
\nc{\BN}{{\mathbb{N}}}
\nc{\BQ}{{\mathbb{Q}}}
\nc{\BP}{{\mathbb{P}}}
\nc{\BR}{{\mathbb{R}}}
\nc{\BZ}{{\mathbb{Z}}}
\nc{\BS}{{\mathbb{S}}}

\nc{\CA}{{\mathcal{A}}}
\nc{\CB}{{\mathcal{B}}}
\nc{\CalC}{{\mathcal C}}
\nc{\CalD}{{\mathcal D}}
\nc{\CE}{{\mathcal{E}}}
\nc{\CF}{{\mathcal{F}}}
\nc{\CG}{{\mathcal{G}}}
\nc{\CH}{{\mathcal{H}}}
\nc{\CK}{{\mathcal{K}}}
\nc{\CL}{{\mathcal{L}}}
\nc{\CM}{{\mathcal{M}}}
\nc{\CMM}{{\mathcal{M}^{\operatorname{gen}}_\hbar(-\rho)}}
\nc{\CN}{{\mathcal{N}}}
\nc{\CO}{{\mathcal{O}}}
\nc{\CP}{{\mathcal{P}}}
\nc{\CQ}{{\mathcal{Q}}}
\nc{\CR}{{\mathcal{R}}}
\nc{\CS}{{\mathcal{S}}}
\nc{\CT}{{\mathcal{T}}}
\nc{\CU}{{\mathcal{U}}}
\nc{\CV}{{\mathcal{V}}}
\nc{\CW}{{\mathcal{W}}}
\nc{\CX}{{\mathcal{X}}}
\nc{\CY}{{\mathcal{Y}}}
\nc{\CZ}{{\mathcal{Z}}}

\nc{\gen}{{\operatorname{gen}}}
\nc{\cM}{{\check{\mathcal M}}{}}
\nc{\csM}{{\check{\mathcal A}}{}}
\nc{\obM}{{\overset{\circ}{\mathbf M}}{}}
\nc{\oCA}{{\overset{\circ}{\mathcal A}}{}}
\nc{\obA}{{\overset{\circ}{\mathbf A}}{}}
\nc{\ooM}{{\overset{\circ}{M}}{}}
\nc{\osM}{{\overset{\circ}{\mathsf M}}{}}
\nc{\vM}{{\overset{\bullet}{\mathcal M}}{}}
\nc{\nM}{{\underset{\bullet}{\mathcal M}}{}}
\nc{\obD}{{\overset{\circ}{\mathbf D}}{}}
\nc{\cp}{{\overset{\circ}{\mathbf p}}{}}
\nc{\ofZ}{{\overset{\circ}{\mathfrak Z}}{}}

\nc{\fa}{{\mathfrak{a}}}
\nc{\fb}{{\mathfrak{b}}}
\nc{\fg}{{\mathfrak{g}}}
\nc{\fgl}{{\mathfrak{gl}}}
\nc{\fh}{{\mathfrak{h}}}
\nc{\fj}{{\mathfrak{j}}}
\nc{\fm}{{\mathfrak{m}}}
\nc{\fn}{{\mathfrak{n}}}
\nc{\fu}{{\mathfrak{u}}}
\nc{\fp}{{\mathfrak{p}}}
\nc{\frr}{{\mathfrak{r}}}
\nc{\fs}{{\mathfrak{s}}}
\nc{\ft}{{\mathfrak{t}}}
\nc{\fT}{{\mathfrak{T}}}
\nc{\ofT}{{\overline{\mathfrak T}}}
\nc{\ofS}{{\overline{\mathfrak S}}}
\nc{\fsl}{{\mathfrak{sl}}}
\nc{\hsl}{{\widehat{\mathfrak{sl}}}}
\nc{\hgl}{{\widehat{\mathfrak{gl}}}}
\nc{\hg}{{\widehat{\mathfrak{g}}}}
\nc{\chg}{{\widehat{\mathfrak{g}}}{}^\vee}
\nc{\hn}{{\widehat{\mathfrak{n}}}}
\nc{\chn}{{\widehat{\mathfrak{n}}}{}^\vee}

\nc{\fA}{{\mathfrak{A}}}
\nc{\fB}{{\mathfrak{B}}}
\nc{\fD}{{\mathfrak{D}}}
\nc{\fE}{{\mathfrak{E}}}
\nc{\fF}{{\mathfrak{F}}}
\nc{\fG}{{\mathfrak{G}}}
\nc{\fI}{{\mathfrak{I}}}
\nc{\fJ}{{\mathfrak{J}}}
\nc{\fK}{{\mathfrak{K}}}
\nc{\fL}{{\mathfrak{L}}}
\nc{\fM}{{\mathfrak{M}}}
\nc{\fN}{{\mathfrak{N}}}
\nc{\frP}{{\mathfrak{P}}}
\nc{\fS}{{\mathfrak S}}
\nc{\fU}{{\mathfrak{U}}}
\nc{\fZ}{{\mathfrak{Z}}}

\nc{\bb}{{\mathbf{b}}}
\nc{\bc}{{\mathbf{c}}}
\nc{\be}{{\mathbf{e}}}
\nc{\bj}{{\mathbf{j}}}
\nc{\bn}{{\mathbf{n}}}
\nc{\bp}{{\mathbf{p}}}
\nc{\bq}{{\mathbf{q}}}
\nc{\bv}{{\mathbf{v}}}
\nc{\bx}{{\mathbf{x}}}
\nc{\by}{{\mathbf{y}}}
\nc{\bw}{{\mathbf{w}}}
\nc{\bA}{{\mathbf{A}}}
\nc{\bB}{{\mathbf{B}}}
\nc{\bC}{{\mathbf{C}}}
\nc{\bK}{{\mathbf{K}}}
\nc{\bD}{{\mathbf{D}}}
\nc{\bH}{{\mathbf{H}}}
\nc{\bM}{{\mathbf{M}}}
\nc{\bN}{{\mathbf{N}}}
\nc{\bS}{{\mathbf{S}}}
\nc{\bT}{{\mathbf{T}}}
\nc{\bV}{{\mathbf{V}}}
\nc{\bW}{{\mathbf{W}}}
\nc{\bX}{{\mathbf{X}}}
\nc{\bP}{{\mathbf{P}}}
\nc{\bZ}{{\mathbf{Z}}}

\nc{\sA}{{\mathsf{A}}}
\nc{\sB}{{\mathsf{B}}}
\nc{\sC}{{\mathsf{C}}}
\nc{\sD}{{\mathsf{D}}}
\nc{\sF}{{\mathsf{F}}}
\nc{\sK}{{\mathsf{K}}}
\nc{\sM}{{\mathsf{M}}}
\nc{\sO}{{\mathsf{O}}}
\nc{\sQ}{{\mathsf{Q}}}
\nc{\sP}{{\mathsf{P}}}
\nc{\sV}{{\mathsf{V}}}
\nc{\sW}{{\mathsf{W}}}
\nc{\sZ}{{\mathsf{Z}}}
\nc{\sfp}{{\mathsf{p}}}
\nc{\sr}{{\mathsf{r}}}
\nc{\sfb}{{\mathsf{b}}}
\nc{\sfc}{{\mathsf{c}}}
\nc{\sd}{{\mathsf{d}}}
\nc{\sg}{{\mathsf{g}}}
\nc{\sfl}{{\mathsf{l}}}

\nc{\BK}{{\bar{K}}}

\nc{\tA}{{\widetilde{\mathbf{A}}}}
\nc{\tB}{{\widetilde{\mathcal{B}}}}
\nc{\tG}{{\widetilde{G}}}
\nc{\TM}{{\widetilde{\mathbb{M}}}{}}
\nc{\tO}{{\widetilde{\mathsf{O}}}{}}
\nc{\tU}{{\widetilde{\mathfrak{U}}}{}}
\nc{\TZ}{{\tilde{Z}}}
\nc{\tZ}{\widetilde{Z}{}}
\nc{\tx}{{\tilde{x}}}
\nc{\tbv}{{\tilde{\bv}}}
\nc{\tfP}{{\widetilde{\mathfrak{P}}}{}}
\nc{\tz}{{\tilde{\zeta}}}
\nc{\tmu}{{\tilde{\mu}}}

\nc{\td}{\ddot{\underline{d}}{}}
\nc{\tzeta}{\widetilde{\zeta}{}}
\nc{\hd}{{\widehat{\underline{d}}}}
\nc{\hG}{{\widehat{G}}}
\nc{\hBP}{\widehat{\mathbb P}{}}
\nc{\hQ}{{\widehat{Q}}}
\nc{\hsM}{\widehat{\mathsf M}{}}
\nc{\hfM}{\widehat{\mathfrak M}{}}
\nc{\hCP}{\widehat{\mathcal P}{}}
\nc{\hCR}{\widehat{\mathcal R}{}}
\nc{\hCS}{{\widehat{\mathcal S}}}
\nc{\hfZ}{\widehat{\mathfrak Z}{}}

\nc{\urho}{\underline{\rho}}
\nc{\uB}{\underline{B}}
\nc{\uC}{{\underline{\mathbb{C}}}}
\nc{\ui}{\underline{i}}
\nc{\ofP}{{\overline{\mathfrak{P}}}}

\nc{\hrho}{{\hat{\rho}}}

\nc{\unl}{\underline}
\nc{\ol}{\overline}
\nc{\one}{{\mathbf{1}}}
\nc{\two}{{\mathbf{t}}}

\nc{\Tot}{{\mathop{\operatorname{\rm Tot}}}}
\nc{\Hilb}{{\mathop{\operatorname{\rm Hilb}}}}
\nc{\CHom}{{\mathop{\operatorname{{\mathcal{H}}\it om}}}}
\nc{\defi}{{\mathop{\operatorname{\rm def}}}}
\nc{\length}{{\mathop{\operatorname{\rm length}}}}

\nc{\Cliff}{{\mathsf{Cliff}}}
\nc{\Fl}{{\mathsf{Fl}}}
\nc{\Fib}{{\mathsf{Fib}}}
\nc{\Coh}{{\mathsf{Coh}}}
\nc{\FCoh}{{\mathsf{FCoh}}}

\nc{\reg}{{\text{\rm reg}}}

\nc{\cplus}{{\mathbf{C}_+}}
\nc{\cminus}{{\mathbf{C}_-}}
\nc{\cthree}{{\mathbf{C}_*}}
\nc{\Qbar}{{\bar{Q}}}
    
\nc{\bh}{{\bar{h}}}
\nc{\bOmega}{{\overline{\Omega}}}
\nc\tGr{\widetilde{\Gr}}

\nc{\seq}[1]{\stackrel{#1}{\sim}}
\nc\ogu{\overline{G/U}}
\nc\chlam{\check{\lam}}

\nc\St{\operatorname{St}}

\nc\uS{\underline{S}}
\nc\QM{\mathcal{QM}}

\begin{document}
\title{Primitive forms for Gepner singularities}
\author{Andrei Ionov}

\begin{abstract}

We provide a construction of Saito primitive forms for Gepner singularity by studying the relation between Saito primitive forms for Gepner singularities and primitive forms for singularities of the form $F_{k,n}=\sum_{i=1}^n x_i^k$ invariant under the natural $S_n$-action.

\end{abstract}

\maketitle

\sec{}{Introduction}{} The {\itshape Gepner singularity} $G_{k,n}$ is the quotient of the singularity of  $F_{k,n}=\sum_{i=1}^n x_i^k$ at the origin by the action of $S_n$ by the permutation of coordinates. 

The interest in the Gepner singularities appeared after the isomorphism of the chiral ring of a $SU(n+1)_{k-n-1}/(SU(n)_{k-n}\x U(1))$ Kazama--Suzuki model, the Milnor ring of the Gepner singularity $G_{k,n}$ and the cohomology ring of the Grassmannian $\Gr(n,k)$ were established  in \cite{G}. The further explorations of the relation between the Gepner singularities and topological conformal field theories (TCFTs) continued in \cite{Z}, \cite{G-ZV}.

All three sides of the isomorphism of \cite{G} admit the natural deformations equipped with a structure of Frobenius manifold: deformations by Witten's descent (\cite{DVV}) for chiral rings of TCFTs, Saito structure for Milnor rings of singularity (\cite{S}) and quantum cohomology for cohomology rings. It appears that the Frobenius structure of quantum cohomology of Grassmannian is not isomorphic to the other two structures even in the simple cases. However, in \cite{DVV} it was proved that the Saito structure for the singularity $z^{k+1}$ is isomorphic to the Frobenius manifold for the $SU(2)_k/U(1)$ Kazama--Suzuki model (also known as minimal models). This leads to a natural conjecture of relation between a Saito structure for Gepner singularity and the Witten's descent deformations of chiral ring of Kazama--Suzuki model formulated in \cite{BGK}. More precisely, there should be a certain Saito primitive form providing Frobenius manifold isomorphic to the one coming from the Witten's descent deformations. Further study of primitive forms for Gepner singularities and corresponding Frobenius structures continued in \cite{BB1}, \cite{BB2}, \cite{BS}, \cite{S0}. 

The notion of a primitive form was introduced in \cite{S} in the setting of versal deformations of a singularity. A choice of a primitive form endows the space of versal deformations of a singularity with a structure of a Frobenius manifold. The key existence theorem of primitive forms for general singularity was proved in \cite{MS}. In \cite{LLS} the dimension of the moduli space of primitive forms for a given singularity was computed. The dimension grows fast as the singularity becomes more complicated. There are only a few examples of explicit constructions of primitive forms.

In the current paper we explore the relation between primitive forms for the Gepner singularity $G_{k,n}$ and for the singularity $F_{k,n}$. More precisely, we use the construction of \cite{CKS} of a Frobenius manifold from a data of a Frobenius manifold with a finite group action to construct primitive forms for $G_{k,n}$ starting from the primitive form for $F_{k,n}$ invariant with respect to natural $S_n$-action. 
It is a natural singularity theory analogue 
 of one of the cases of the abelian/nonabelian correspondence in  \cite{CKS} relating the quantum cohomology of the  Grassmannian $\Gr(n,k)$ and the product of projective spaces $(\BP^{k-1})^{\x n}$.

Remarkably, the relation we study should also impose a relation between the $SU(n+1)_{k-n-1}/(SU(n)_{k-n}\x U(1))$ Kazama--Suzuki model and the tensor product of $n$ copies of the minimal $SU(2)_k/U(1)$ model, which is to be investigated.

One of the interesting prospects for the work would be the generalization of the results of the paper to the case of a singularity invariant under the action of a complex reflection group and the corresponding quotient. Another interesting question is to understand the relations with the equivariant singularity theory. 

\refs{s}---\refs{pf} consist of definitions and preliminary facts. In \refs{W} we introduce the construction of \cite{CKS} for Frobenius manifold with a finite group action. In \refs{m}, \refs{p} we construct primitive forms for Gepner singularities. In \refs{e} we compare our construction with previously known constructions of primitive forms in the known cases.

\ssec{}{Acknowledgements} The author is grateful to A. Belavin, K. Hori, Y. Kononov, T. Milanov, S. Natanzon and K. Saito for the useful discussions. The author is grateful to Kavli IPMU and especially K. Saito for hospitality during the visits in 2015 and 2016. The author is grateful to D. Kubrak and S. Natanzon for reading the original manuscript of the paper. 

The study has been funded
by the Russian Academic Excellence Project '5-100'. The author was supported  in part by
Dobrushin stipend and grant RFBR 15-01-09242.

\sec{s}{Preliminaries on singularity theory}

Let $\bfz=(z_0,\ldots, z_n)$, let $\CC\{\bfz\}$ be the ring of germs of holomorphic functions in $\bfz$ at the origin and let $f=f(\bfz)\in \CC\{\bfz\}$ be a function with an isolated singularity at the origin. The Milnor ring of $f$ is defined to be a quotient $J_f:=\CC\{\bfz\}/I_f,$ where $I_f:=(\frac{\pa f}{\pa z_0},\ldots,\frac{\pa f}{\pa z_n})$ is a Milnor ideal. Under the isolated singularity assumption $\mu:=\dim J_f$ is finite and is called the Milnor number.

Let $V$ be a vector space with coordinates $\bft=(t_1,\ldots,t_m)$. The germ of a holomorphic function at the origin $F=F(\bfz,\bft)\in\BC\{\bfz, \bft\}$ is said to be a deformation of $f$ if $F(\bfz,0)=f(\bfz)$ for all $\bfz$ as germs at the origin.

Let $(V,F)$ be a deformation of $f$ and let $\CT_V$ be a $\CC\{\bft\}$-module of germs of holomorphic vector fields at the origin of $V$. Then there is a well defined Kodaira--Spencer map of $\CC\{\bft\}$-modules: $$\mathrm{KS}\colon\CT_V\to \BC\{\bfz, \bft\}/\Big(\frac{\pa F}{\pa z_0},\ldots,\frac{\pa F}{\pa z_n}\Big),$$
defined as follows. For a germ of vector field $\xi$ pick its lift $\tilde{\xi}$ to a germ on $V\x \BC^{n+1}$ and put $\mathrm{KS}(\xi)$ to be equal to the image of a derivative $\tilde{\xi}(F)$.

The deformation $(V,F_1)$ is said to be induced from $(U,F_2)$ with respect to holomorphic map $h\colon V\to U$ mapping the origin in $V$ to the origin in $U$ if $F_1(\bfz,\bft)=F_2(\bfz,h(\bft))$. The deformation $(V_f, \tilf)$ is said to be {\itshape versal} if every deformation is induced from it and it is of minimal dimension among such.

The versal deformation always exists and is unique up to isomorphism, moreover $V_f\simeq J_f$, in particular $\dim V_f=\mu$. The Kodaira--Spencer map is an isomorphism if and only if the deformation is versal.

For proofs and further details we refer to \cite{AG_ZV}.

\ssec{}{Gepner singularities} Let us fix two positive integers $n$ and $k$ and let $F_{k,n}=\sum_{i=1}^n x_i^k$ be a polynomial in variables $\bfx=(x_1,\ldots,x_n)$. Let us define $\bfy=(y_1,\ldots,y_n)$ by $y_i=\sig_i(\bfx)$, where $\sig_i$ is the $i$-th elementary symmetric function in $n$ variables, so that $$1+\sum_{i=1}^n y_i T^i=\prod_{i=1}^n (1+x_iT).$$ Since the polynomial $F_{k,n}$ is symmetric, there is a polynomial $G_{k,n}$, such that $F_{k,n}(x_1,\ldots,x_n)=G_{k,n}(y_1,\ldots,y_n)$. 
It is not hard to verify that

\prop{}{}(\cite{G},\cite{S0})
If $k>n$ then $G_{k,n}$ has an isolated singularity at the origin $\bfy=0$. In this case the Milnor number of ${G_{k,n}}$ is equal to $\binom{k-1}{n}$.
\eprop

In what follows we will assume that $k>n$. We will call this singularity a {\itshape Gepner singularity}.

\sec{F}{Frobenius manifolds} 

Let $V$ be a finite dimensional vector space over a field $\bfk$ of characteristic $0$ and let $V^\vee$ be its dual. Then let us put $M=Spf(\bfk[[V^\vee]])$ to be the formal completion of $V$ at the origin, so that the functions on $M$ are formal series in $V^\vee$. We denote by $T_M$ its tangent sheaf which is canonically isomorphic to $V\ten\CO_M$.

\defe{}{}
{\itshape The formal Frobenius manifold} on $M $
 is the collection of data: $(\bullet, g, e, \CE)$ 
 , where

1) $g$ is a $\CO_M$-linear nondegenerate pairing on $T_M$ such that the corresponding connection $\nabla$ is flat;

2) $\bullet$ is $\CO_M$-linear, associative, commutative product on $T_M$, such that $\nabla c$ is symmetric where $c$ is the tensor defined as $c(u,v,w)=g(u\bullet v,w)$;

3) $e$ is a formal vector field on $M$, which is the identity for $\bullet$ and such that $\nabla e=0$;

4) $\CE$ is a formal vector field  on $M$, which is called an Euler vector field, and satisfies
$$\nabla\nabla\CE=0,\enskip\CL_\CE g=Dg, \enskip\CL_\CE (\bullet)=\bullet, \enskip\CL_\CE(e)=-e,$$
where $\CL$ denote the Lie derivative and $D\in\bfk$ is a constant.
\edefe

For further discussions of the notion we refer to \cite{D}, \cite{M1}, \cite{Sa}.

\sec{pf}{Saito structures and primitive forms}

We introduce the notions of {\itshape pre-Saito structure} and {\itshape a primitive section} following \cite{Sa}.

\defe{}{}
 We will call {\itshape pre-Saito structure} the following data: $(M, E, g, \nabla, \Phi, R_0, R_\infty) $, where 

1) $M$ is a formal completion at the origin of a finite-dimensional vector space $V$ over field $\bfk$ of characteristic $0$;

2) $E$ is a free $\CO_M$-module of finite rank with a flat connection $\nabla$ and a $\CO_M$-bilinear form $g$ flat with respect to $\nabla$;

3) $\Phi, R_0$ and $R_\infty$ are $\CO_M$-linear morphisms $\Phi\colon T_M\ten_{\CO_M}E\to E$ and $R_0, R_\infty \colon E\to E$, satisfying the conditions $$\nabla_{\partial_{t_i}}\Phi_{\partial_{t_j}}=\nabla_{\partial_{t_j}}\Phi_{\partial_{t_i}}, \enskip [\Phi_{\partial_{t_j}},\Phi_{\partial_{t_i}}]=0, \enskip [R_0,\Phi_{\partial_{t_i}}]=0,$$
$$\nabla(R_\infty)=0, \enskip \Phi_{\partial_{t_i}}+\nabla_{\partial_{t_i}}R_0=[\Phi_{\partial_{t_i}},R_\infty],$$ $$\Phi_{\partial_{t_i}}^*=\Phi_{\partial_{t_i}}, \enskip R_0^*=R_0, \enskip R^*_\infty+R_\infty=-w\Id, $$
where $w\in\BZ$ is a fixed integer called the weight, $*$ stands for $g$-adjoint, $\{t_i\}$ are the coordinates on $M$ induced by a basis of   $V$ and $\Phi_\xi\colon E\to E$ is the map obtained by the substitution of the section $\xi$ of $T_M$ into $\Phi$.

\edefe

\defe{}{} A section $\ome\in\Gam(M,E)$ is called {\itshape a homogeneous primitive section} if 

1) it is flat: $\nabla(\ome)=0$, 

2) the morphism $\phi_\ome\colon T_M\to E$ given by $\xi\mapsto \Phi_\xi(\ome)$ is an isomorphism and 

3) $R_\infty \ome=q\ome$ for some $q\in\bfk$. 
\edefe

Given a pre-Saito structure $(M, E, g, \nabla, \Phi, R_0, R_\infty) $ and a homogeneous primitive section $\ome$ for it, one constructs a structure of Frobenius manifold on $M$ by taking $^\ome\nabla:=\phi^{-1}_\ome\nabla\phi_\ome$, $\xi\bullet\eta:=-\Phi_\xi(\phi_\ome(\eta))$, $e:=\phi_\ome^{-1}(\ome)$, $\CE:=\phi^{-1}_\ome(R_0(\ome))$, $^\ome g(\xi,\eta):=g(\phi_\ome(\xi),\phi_\ome(\eta))$.

\ssec{}{Primitive forms for an isolated singularity (\cite{S}, \cite{ST})} 

Let us return to the notations of \refs{s}. Let $M_f$ be a formal completion of $V_f$ at the origin. Note that the Kodaira--Spencer isomorphism endows $T_{M_f}$ with a $\CO_M$-bilinear product $\bullet$. We also define vector fields $\CE:=\mathrm{KS}^{-1}(\tilf)$ and $e:=\mathrm{KS}^{-1}(1)$.

Consider a $\BC\{\bft\}$-module consisting of germs of forms of top degree in $\bfz$-variables $\varphi(\bfz,\bft)dz_0\wedge\ldots\wedge dz_n$ modulo the image of the wedge multiplication by $1$-form $d\tilf$. After passing to the formal completion at the origin we obtain a vector bundle $\Ome_\tilf$ on $M_f$. It possesses  the bilinear residue pairing 
$$
g(\ome_1, \ome_2):=\Res\left[
  \begin{array}{cccc}
   \varphi_1\varphi_2dz_0\wedge\ldots\wedge dz_n \\
  \frac{\pa \tilf}{\pa z_0},\ldots,\frac{\pa \tilf}{\pa z_n}
  \end{array}
\right]\in\CO_{M_f},
$$
for $\ome_i=\varphi_idz_0\wedge\ldots\wedge dz_n$.

The multiplication of a form by a function together with the Kodaira--Spencer isomorphism provides a bilinear map $\Phi\colon T_{M_f}\ten_{\CO_{M_f}}\Ome_\tilf\to \Ome_\tilf$. Then, as above, a form $\ome\in\Gam(M_f,\Ome_\tilf)$ defines a map $\phi_\ome\colon T_{M_f}\to \Ome_\tilf$ given by $\xi\mapsto \Phi_\xi(\ome)$. We then put $^\ome\! g(\xi,\eta):=g(\phi_\ome(\xi),\phi_\ome(\eta))$.

We will call such $\ome$  {\itshape a primitive form} if  $(\bullet,  ^\ome\!\!g, e, \CE)$ provides a Frobenius manifold structure on $M_f$.  Primitive forms always exist (\cite{MS}) but, in general, are not unique.

If $f=z^{k+1}$ then the class $dz\in\Gam(M_{z^{k+1}},\Ome_{\widetilde {z^{k+1}}})$ is the unique (up to scalar multiplication) primitive form. We will call the corresponding Frobenius manifold $\CA_k$.

\sec{W}{Frobenius manifold with finite group action}
Let $(M,\bullet,g,e,\CE)$ be a Frobenius manifold and let $W$ be a finite group acting on $M$ by automorphisms in a way compatible with the Frobenius structure.

Let us consider the fixed point set of the $W$-action $M^W$. Then $M^W$ is a smooth formal subscheme of $M$,  $W$ acts $\CO_{M^W}$-linearly on $T_M|_{M^W}$ and $T_{M^W}=(T_M|_{M^W})^W$.

Let us fix a non-trivial character $\sgn\colon W\to\pm1$ and consider the corresponding antisymmetrization morphism $a\colon T_M|_{M^W}\to T_M|_{M^W}$ given by $a(\xi)=\sum_{w\in W} \sgn(w)w(\xi)$. We denote its image by $E$.
It is a locally free $\CO_M$-module. Note that 
we have a $g$-orthogonal direct sum decomposition $T_M|_{M^W}=\ker a\oplus E$ and the restriction of $g$ to $E$ is nondegenerate. We denote by $\nabla$ the restriction of the connection on $T_M|_{M^W}$ to $E$.

There is a natural $\CO_{M^W}$-linear multiplication $\Phi\colon (T_M|_{M^W})^W\ten E\to E$ coming from multiplication on $T_M$ and the operator $R_0:=\CE\bullet$. We also put $R_\infty:=\nabla \CE-\Id$. It is easy to check that $R_0$ and $R_\infty$ preserve $E$ and

\lem{}{}(\cite{CKS}, Lemma 2.3.1)
The above $(M^W,E,g,\nabla,\Phi,R_0, R_\infty)$ is a pre-Saito structure.
\elem

We then have

\prop{CKS}{}(\cite{CKS}, Proposition 2.3.2)
Suppose there is a $\nabla$-horizontal $R_\infty$-eigensection $\ome\in\Gam(M^W,E)$ such that the morphism $\phi_\ome\colon T_{M^W}\to E$ given by $\xi\mapsto\xi\bullet\ome$ is surjective. Then every smooth formal subscheme $N\subset M^W$ such that the restriction of the above morphism $T_N\to E$ is an isomorphism has a natural Frobenius structure.
\eprop

\sec{m}{Main Construction} 
Consider the tensor product of Frobenius manifolds $M_{k,n}:=\CA_{k-1}^{\ten n}$ (see \cite{M1} for the definition). Note that the underlying space is naturally identified with 
$M_{F_{k,n}}$. 

\rem{M} It follows from \cite{M}, Theorem 3.2.3, that 
there is a primitive form for $F_{k,n}$, which provides the Frobenius manifold $M_{F_{k,n}}$. 
\erem

The Frobenius manifold $M_{k,n}$ naturally comes with an action of $W=S_n$ by permutation of the factors. We now apply the construction of \refs{W} to it.

\prop{p} 
There is a section $\ome$ of $E$ satisfying the conditions of \refp{CKS}. 
\eprop
\prf

Since we work locally at the origin it follows from Nakayama lemma as in \cite{Sa} Remark VII.3.7 that it is sufficient to construct $\ome$ at the origin and use the flat connection $\nabla$ to translate it.

Consider an antisymmetric polynomial $w_n:=\prod_{1\le i<j\le n} (x_i-x_j)$ as an element of the Milnor ring $J_{F_{k,n}}$. Note that, since $k>n$ we have $w_n\ne0$. Moreover, it is a homogeneous element. It can be viewed as an element of $E_0$, the fibre of $E$ at the origin. At the origin the map $\phi_\ome|_0\colon J^W_{F_{k,n}}\to E_0$ is obviously surjective. Then the statement follows.

\epr

Let us now choose a subscheme $N\subset M_{k,n}^W$ appropriate for the application of \refp{CKS}. We will start with the following

\prop{s}
There is a short exact sequence:
\eq{ex}0\to \ker(w_n\cdot)\to J_{F_{k,n}}^W\to J_{G_{k,n}}\to0,\eeq
where $J_{F_{k,n}}^W$ is the subring of $W$-invariants in the Milnor ring $J_{F_{k,n}}$ and $\ker (w_n\cdot)$ is a kernel in $J_{F_{k,n}}^W$ of multiplication by $w_n\in J_{F_{k,n}}$ 
$$ w_n\cdot\colon J_{F_{k,n}}^W\to J_{F_{k,n}}$$
(cf. (3.1.2) in \cite{CKS}).  
\eprop
\prf
The first arrow is the natural embedding. Let us construct the second arrow. Let $I_{F_{k,n}}\subset \BC\{\bfx\}$ be the Milnor ideal of $F_{k,n}$, let $I_{F_{k,n}}^W\subset \BC\{\bfy\}$ be its $W$-invariant part and let $I_{G_{k,n}}\subset \BC\{\bfy\}$ be the Milnor ideal of $G_{k,n}$. To obtain a surjective map $J_{F_{k,n}}^W\to J_{G_{k,n}}$ it is sufficient to prove that $I_{F_{k,n}}^W\subset I_{G_{k,n}}$. Note that by the chain rule we have: $$\frac{\partial F_{k,n}}{\partial x_i }=\sum_{j=1}^{n}\frac{\partial \sig_j(\bfx)}{\partial x_i }\frac{\partial G_{k,n}}{\partial y_j}(\sig(\bfx)).$$
Therefore, we have $I_{F_{k,n}}\subset I_{G_{k,n}}\BC\{\bfx\}$. Also, we, obviously have $I_{F_{k,n}}^W\BC[\bfx]\subset I_{F_{k,n}}$. Thus, $I_{F_{k,n}}^W\BC\{\bfx\}\subset  I_{G_{k,n}}\BC\{\bfx\}$ and $I_{F_{k,n}}^W\subset I_{G_{k,n}}$.

It remains to check the exactness in the middle term of the sequence. To show that the composition of the two arrows is zero it is sufficient to show that $w_nI_{G_{k,n}}\BC\{\bfx\}\subset I_{F_{k,n}}$. But, the determinant of the Jacobi matrix with the entries $\frac{\partial \sig_j(\bfx)}{\partial x_i }$ is equal to $w_n$. Therefore, we have $$w_n\frac{\partial G_{k,n}}{\partial y_j}(\sig(\bfx))=\sum_{j=1}^{n} a_{ij} \frac{\partial F_{k,n}}{\partial x_i },$$ for some polynomials $a_{ij}\in\BC\{\bfx\}$ (minors of the Jacobi matrix) and the embedding $w_nI_{G_{k,n}}\BC\{\bfx\}\subset I_{F_{k,n}}$ follows.

Finally, the dimension of the cokernel of the first arrow is equal to the dimension of the antiinvariants of $W$ in $J_{F_{k,n}}$, which is equal to the dimension of the space of the antisymmetric polynomials in $\bfx$ modulo $x_i^{k-1}$. And this is equal to the number of monomials of the form $x_1^{i_1}x_2^{i_2}\ldots x_n^{i_n}$ with $k-1>i_1>i_2>\ldots>i_n\ge 0$. By simple combinatorial calculation this is $\binom{k-1}{n}$. It is the same as the dimension of $J_{G_{k,n}}$ and the proposition follows.

\epr

Let us now choose a splitting of the short exact sequence \refe{ex} as a sequence of vector spaces: $\iot\colon J_{G_{k,n}}\to J_{F_{k,n}}^W$. 
The map $\iot$ naturally gives us a formal subscheme $N\simeq M_{G_{k,n}}\subset M_{k,n}^W$. 

Note that the same argument as in proof of \refp{s} together with the Nakayama lemma 
provides a short exact sequence 
\eq{ex1}0\to \ker(\ome\cdot)\to (T_{M_{k,n}}|_{M_{k,n}^W})^W\to T_N\to0\eeq of bundles over $N$. Therefore, we have an isomorphism $\phi_\ome\colon T_N\iso E=\mathrm{im}(\ome\cdot)$. This implies

\lem{sub}
The above $N$ satisfies the conditions of \refp{CKS}.

\elem

We obtain

\th{m}
There is a Frobenius structure on $M_{G_{k,n}}$ depending on the choice of the splitting of the short exact sequence \refe{ex}. Moreover, the data of $(\bullet, e, \CE)$ (with the metric omitted) for this Frobenius manifold and a Frobenius manifold provided by a primitive form for $G_{k,n}$ coincide.
\eth
\prf
The first half follows from \refp{CKS}, \refp{p} and \refl{sub}. Let us prove the second part.

By construction, the second maps of \refe{ex} and \refe{ex1} are algebra homomorphisms with the multiplication on $T_N$ being the multiplication on a Milnor ring of the deformed polynomial $\tilG_{k,n}$. Then the map $\Phi\colon T_N\ten E\to E$ provides a free rank $1$ module structure over $T_N$ on $E$ and the multiplication on $T_N$ provided by \refp{CKS} is the same as the multiplication on $T_N$ described above. Thus, the multiplication $\bullet$ and the identity vector field $e$ coincide for these two Frobenius manifolds.

Since the image of $\tilF_{k,n}$ in $T_N$ under the second map in \refe{ex1} is equal to the image of $\tilG_{k,n}$ in $T_N$ the Euler fields coincide.
\epr

\conj{}
The above Frobenius structure does not depend on the choice of the splitting $\iot$.
\econj

\sec{p}{Primitive forms for Gepner singularities}

Let us fix a splitting in \reft{m}. In this section we prove the following result

\th{prim}
There is a primitive form $\zet$ for $G_{k,n}$ such that the induced Frobenius structure on $M_{G_{k,n}}$ is isomorphic to the Frobenius structure of \reft{m}. 
\eth
\prf
We only need to provide the compatibility of the metrics.

By \refr{M} there is the primitive form for $F_{k,n}$ providing the Frobenius structure of $M_{k,n}$. Let us denote by $\eta=\varphi(\bfx, \bft) dx_1\ldots dx_n \in \Gam(N,\Ome_{\tilF_{k,n}})$ its restriction to $N$. The $S_n$-equivariance
 of $M_{k,n}$ implies that $\varphi$ is symmetric in $x_i$.

\lem{j}
There is an isomorphism $j\colon\Ome_{\tilG_{k,n}}\iso \enskip \Ome^W_{\tilF_{k,n}}|_N$ preserving the residue pairings. 
\elem
\prf
We define $j$ to be a morphism induced by the change of variables $\{\bfx\}\mapsto\{\bfy\}$. More precisely, let $\psi(\bfy,\bft)dy_1\ldots dy_n$ be an element of $\Ome_{\tilG_{k,n}}$ then $j(\psi(\bfy,\bft)dy_1\ldots dy_n)=\psi(\bfx,\bft)w_n dx_1\ldots dx_n\in\Ome^W_{\tilF_{k,n}}|_N$. The map is well defined since $d\tilF_{n,k}$ equals $d\tilG_{n,k}$ on $N$ after the change of variables. 
The map is an isomorphism, since sections  of $\Ome^W_{\tilF_{k,n}}|_N$ are exactly sections $\theta(\bfx,\bft)w_n dx_1\ldots dx_n$ of $\Ome_{\tilF_{k,n}}|_N$ with $\theta(\bfx,\bft)$ symmetric in $x_i$.
Verification of compatibility with the pairing is straightforward.

\epr
Consider now the diagram 

$$
\xymatrix{ T_{M_{k,n}}|_N\ar[r]^{\phi_\eta}_\sim & \Ome_{\tilF_{k,n}}|_N \\
E|_N\ar[r]_\sim^{\phi_\eta} \ar@{^{(}->}[u] &  \Ome^W_{\tilF_{k,n}}|_N\ar@{^{(}->}[u] \\
T_N\ar[u]^{\phi_\ome}_\sim \ar@{-->}[r] & \Ome_{\tilG_{k,n}}\ar[u]^{j}_{\sim}. 
}$$ 
We define the lower arrow to be the composition $j^{-1}\circ\phi_\eta\circ\phi_\ome$. It is easy to see that this map is $\phi_\zet$ for $\zet=\frac{\ome}{w_n} \varphi(\bfy,\bft) dy_1\ldots dy_n.$ Moreover, it follows from the diagram and \refl{j} that the Frobenius structure of \reft{m} is induced by $\zet$. This implies the theorem.
\epr

\rem{} One can construct the corresponding element of the filtered de Rham complex (see \cite{S}, \cite{ST} for the definition) in a similar way.
\erem

\sec{e}{Examples of primitive forms}

Let us look at what our results provide in the simplest examples. Let us first list the simplest singularities among Gepner singularities. It follows immediately from \cite{S0} that

\th{} a) The singularities $G_{k,1}, G_{k+1, k}, G_{k+2, k}, G_{5,2}$ are the only simple (\cite{AG_ZV}) singularities among all $G_{k,n}$. In these cases they have respectively the types $A_{k-1}, A_1, A_{k+1}, D_6$.

b) The singularities $G_{6,2}$ and $G_{6,3}$ are the only unimodal (\cite{AG_ZV}) singularities among all $G_{k,n}$. In these cases they are simple elliptic singularities of the type $\tilE_8$, i.e. they are equivalent up to stabilization to the hypersurface singularity given by $x^2+y^3+z^6+\sig y^2z^2=0$ for some fixed values of parameter $\sig$, such that $4\sig^3+27\ne 0$.
\eth

It follows from \cite{LLS} that simple singularities are the only singularities for which there is a unique primitive form (up to multiplicative constant) and unimodal singularities are the only singularities for which the moduli space of primitive forms (up to multiplicative constant) is one-dimensional. It, therefore,  follows that if $G_{k,n}$ is a simple singularity then the primitive form provided by \reft{m} and \reft{prim} is the unique primitive form for this singularity and if $G_{k,n}$ is a unimodal singularity this primitive form is a point in the one-dimensional space of all primitive forms.

The primitive forms for the simple elliptic singularities are constructed in \cite{S}. For the singularity of the type $\tilE_8$ the construction goes as follows. To the family of marginal deformations $x^2+y^3+z^6+\sig y^2z^2=0$ one associates the family of elliptic curves $E_\sig$. Then for a choice of a cycle $A\in H_1(E_\sig,\BC)$ one defines a function $\pi_A(\sig)$ as a certain period integral over cycle $A$. This function satisfies a hypergeometric Picard--Fuchs
equation. Now the primitive form is given by $\zeta=\zeta(\sig)=\frac{dx\wedge dy\wedge dz}{\pi_A(\sig)}$. We refer to \cite{S} for more details. 

It is now natural to ask 

\noindent
{\bfseries Question 8.2.} For unimodal Gepner singularities $G_{6,2}$ and $G_{6,3}$, what are the relations between the choice of a splitting in \reft{m} and the choice of a cycles $A$ in \cite{S}?

\bigskip

\footnotesize{
{\bf A.I.}: National Research University
Higher School of Economics, Russian Federation, Department of Mathematics, 6 Usacheva st, Moscow 119048;\\
Massachusetts Institute of Technology, Department of Mathematics, 77 Massachusetts
Ave, Cambridge, MA 02139 USA;\\
{\tt aionov@mit.edu}}

\end{document}